\documentclass[11pt,twoside]{article}

\usepackage{amsfonts}
\usepackage[mathscr]{euscript}

\newcommand{\lbl}[1]{\label{#1}}

\newtheorem{theo}{Theorem}[section]
\newtheorem{prop}{Proposition}[section]
\newtheorem{lem}{Lemma}[section]
\newtheorem{remark}{Remark}[section]

\newcommand{\be}{\begin{equation}}
\newcommand{\ee}{\end{equation}}
\newcommand\bes{\begin{eqnarray}} \newcommand\ees{\end{eqnarray}}
\newcommand{\bess}{\begin{eqnarray*}}
\newcommand{\eess}{\end{eqnarray*}}
\newcommand\bedd{\bes\left\{\begin{array}{ll}\medskip}
\newcommand\eedd{\end{array}\right.\ees}

\newcommand\ep{\varepsilon}
\newcommand\lk{\left}
\newcommand\rr{\right}
\newcommand\nm{\nonumber}
\newcommand\dd{\displaystyle}
\newcommand\vp{\varepsilon}

\newcommand\ud{\underline}
\newcommand\la{\lambda}

\newcommand\bbb{\big}
\newcommand\ff{, \ \ \forall \ }
\begin{document}
  \pagestyle{myheadings}
\thispagestyle{empty}

\title{Free boundary problems for a Lotka-Volterra competition system\thanks{This work was
supported by NSFC Grants 11071049 and 11371113}}

\author{{\Large Mingxin Wang}\thanks{Natural Science Research Center, Harbin Institute of Technology, Harbin 150080, PR China.}, \ \ \
{\Large Jingfu Zhao}
\thanks{Natural Science Research Center, Harbin Institute of Technology, Harbin 150080, PR China, and
Department of Mathematics, Shaanxi University of Technology, Hanzhong
723000, PR China. }}
\date{2013/01/22}
\maketitle

\begin{quote}
\noindent{\bf Abstract.} In this paper we investigate two free boundary problems for a Lotka-Volterra type competition model in one space dimension. The main objective is to
understand the asymptotic behavior of the two competing species spreading
via a free boundary. We prove a spreading-vanishing dichotomy, namely the two species
either successfully spread to the right-half-space as time $t$ goes to infinity and survive
in the new environment, or they fail to establish and die out in the long run. The
long time behavior of the solutions and criteria for spreading and vanishing are also
obtained. This paper is an improvement and extension of J. Guo and C. Wu \cite{GW}.


 \end{quote}

 \section{Introduction}
 \setcounter{equation}{0} {\setlength\arraycolsep{3pt}

We study the evolution of positive solutions $(u,v,s)$ to the following free boundary problems for Lotka-Volterra type competition system
$$\left\{\begin{array}{lll}
 u_t=u_{xx}+u(1-u-kv), &t>0, \ \ 0<x<s(t),\\[1mm]
 v_t=Dv_{xx}+rv(1-v-hu),&t>0, \ \ 0<x<s(t),\\[1mm]
 u_x=v_x=0,\ \ &t>0, \ \ x=0,\\[1mm]
 u=v=0, \ s'(t)=-\mu(u_x+\rho v_x),\ \ &t>0, \ \ x=s(t), \\[1mm]
 u(0,x)=u_0(x), \ \ v(0,x)=v_0(x),& 0\leq x\leq s_0,\\[1mm]
 s(0)=s_0
 \end{array}\right.\eqno{\rm(NFB)}
$$
 and
$$\left\{\begin{array}{llll}
 u_t=u_{xx}+u(1-u-kv), &t>0, \ \ 0<x<s(t),\\[1mm]
 v_t=Dv_{xx}+rv(1-v-hu),&t>0, \ \ 0<x<s(t),\\[1mm]
 u=v=0,\ \ &t>0, \ \ x=0,\\[1mm]
 u=v=0, \ s'(t)=-\mu(u_x+\rho v_x),\ \ &t>0, \ \ x=s(t), \\[1mm]
 u(0,x)=u_0(x), \ \ v(0,x)=v_0(x),& 0\leq x\leq s_0,\\[1mm]
 s(0)=s_0.
 \end{array}\right.\eqno{\rm(DFB)}
$$
In the above two problems, $x=s(t)$ represents the moving boundary, which is to be determined, $k, h, r, D, s_0, \mu$ and $\rho$ are given positive constants. The initial functions $u_0(x),v_0(x)$ satisfy

(NFB1)\, $u_0,\,v_0\in C^2([0,s_0])$, $u_0'(0)=v_0'(0)=u_0(s_0)=v_0(s_0)=0$, $u_0(x), v_0(x)>0$ in  $(0,s_0)$ for the problem (NFB);

(DFB1)\, $u_0,\,v_0\in C^2([0,s_0])$, $u_0(0)=v_0(0)=u_0(s_0)=v_0(s_0)=0$, $u_0(x), v_0(x)>0$ in  $(0,s_0)$ for the problem (DFB).

Problems (NFB) and (DFB) may be viewed as describing the spreading of two new or invasive competing species with population density $(u(t,x),v(t,x))$ over a one dimensional habitat.  In general, both species have a tendency to
emigrate from the boundary to obtain their new habitat, i.e., they will move outward
along the unknown curve (free boundary) as time increases. It is assumed that the expanding speed of the free boundary is proportional to the normalized population gradients at the free boundary, i.e.,
 \[ s'(t)=-\mu\big[u_x(t,s(t))+\rho v_x(t,s(t))\big],\]
which is the well-known Stefan type condition and whose ecological background can refer to \cite{BDK}. Such kind of free boundary conditions have been used in \cite{GW, HMS, Lin, MYS1, MYS2, MYS3}.

Recently, Guo and Wu \cite{GW} studied the problem (NFB) with weak competition case: $0<k,h<1$. Let $s_\infty=\lim\limits_{t\to\infty}s(t)$, $s_*=\frac \pi2\min\bbb\{\sqrt{D/ r},\, 1\bbb\}$,
 \[s^*=\dd\left\{\begin{array}{ll}
 \dd\frac{\pi}{2}\sqrt{\frac D r}\frac 1 {\sqrt{1-h}},\ \ &D<r,\\[4mm]
 \dd\frac\pi 2\frac 1 {\sqrt{1-k}},\ \ &D>r,\\[4mm]
 \dd\frac\pi 2\min\lk\{\frac 1 {\sqrt{1-k}},\frac 1 {\sqrt{1-h}}\rr\},\ \ &D=r,
 \end{array}\right.
 \]
 and
\[
\begin{array}{l}
A=\left\{0<D<r,\ 0<h<1-\frac D r,\ 0<k<1,\ \mu,\rho>0\right\},\\[2mm]
B=\left\{0<r<D,\ 0<k<1-\frac r D,\ 0<h<1,\ \mu,\rho>0\right\}.
 \end{array}
\]
They obtained the following results:

(a) If $s_\infty<s_*$, then the two species vanish eventually, i.e.,
  \[\lim_{t\to\infty}\|u(t,\cdot)\|_{C([0,s(t)])}= \lim_{t\to\infty}\|v(t,\cdot)\|_{C([0,s(t)])}=0;\]
if $s_\infty>s^*$, then the two species spreading successfully, i.e., $\dd\lim_{t\to\infty} u(t,\cdot)>0$, $\dd\lim_{t\to\infty} v(t,\cdot)>0$.

\vskip 4pt (b) If $(D,h,k,r,\mu,\rho)\in A\bigcup B$, then either $s_\infty<s_*$ (and so the two species vanish eventually), or the two species spreading successfully.

(c) If the two species spreading successfully, then
  \[\lim_{t\to\infty}(u,v)(t,x)=\left(\frac{1-k}{1-hk}, \frac{1-h}{1-hk}\right)\]
uniformly in any compact subset of $[0,\infty)$.

The main purpose of this paper is to improve and extend the above results obtained in \cite{GW}. We shall remove the restriction $0<k,h<1$ and give a complete description for the spreading-vanishing dichotomy, long time behavior of $(u,v)$ and sharp criteria for spreading and vanishing.

For the global existence and uniqueness of the solution, from the proof of theorem 1 in \cite{GW}, we can easily see that any one of both problems (NFB) and (DFB) admits a unique global solution $(u,v,s)\in C^{1,2}(\Omega)\times C^{1,2}(\Omega)\times C^1([0,\infty))$ where $\Omega=\{(t,x):\ t>0,\ 0\leq x\leq s(t)\}$. Moreover, the following inequalities hold
  \begin{equation}
0<u(t,x)\leq M,\ 0<v(t,x)\leq M,\ 0<s'(t)\leq\mu M\lbl{1.1}
  \end{equation}
for $t>0,\ 0<x<s(t)$, where the positive constant $M$ depending only on $D$, $r$, $\rho$, $\|u_0\|_{L^\infty}$, $\|v_0\|_{L^\infty}$, $\min_{x\in[0,s_0]}u_0'(x)$ and $\min_{x\in[0,s_0]}v_0'(x)$.

\vskip 4pt Recently, Wang and Zhao \cite{WZ} studied a free boundary problem for a predator-prey model with double free boundaries in one dimension, in which the prey lives in the whole space but the predator lives in a bounded area at the initial state.
Later on, Zhao and Wang \cite{ZW} extended the results to the case of higher space dimension with radially symmetric parameters, and Wang \cite{W} dealt with the case that both predator and prey live in a bounded area at the initial state. They established the spreading-vanishing dichotomy, long time behavior of the solution and sharp criteria for spreading and vanishing.

In the absence of $v$, the problems (NFB) and (DFB) are reduced to the one phase Stefan problems which were studied by Du and Lin in \cite{DLin}, Kaneko and Yamadae in \cite{KY} respectively. The well-known Stefan condition has been used in the modeling of a number of applied problems. For example, it was used to describe the melting of ice in contact with water \cite{RU}, the modeling of oxygen in the muscle \cite{CR}, the wound healing \cite{ChenA}, the tumor growth \cite{ChenF}, and the spreading of species \cite{DLin, DG, HMS, Lin}. There is a vast literature on the Stefan problems, and some important theoretical advances can be found in \cite{CS, CR} and the references therein.

Some similar free boundary problems have been used in two-species models  over a bounded spatial interval in several earlier papers; please refer to, for example,  \cite{HIMN,HMS,Lin,MYS1,MYS2,MYS3}. For the study of free boundary problems for other type biological models, we refer to, for instance \cite{ChenF,DG1, DGP, DLou, PZ} and references cited therein.

The organization of this paper is as follows. Section $2$ is devoted to the long time behavior of $(u,v)$.  From those results we can also get a spreading-vanishing dichotomy. In Section 3 we shall give the sharp criteria for spreading and vanishing. The last section is a brief discussion.

\section{Long time behavior of $(u,v)$}
\setcounter{equation}{0}

It follows from (\ref{1.1}) that $x=s(t)$ is
monotonic increasing. Therefore, there exists $s_\infty\in(0,\infty]$ such that $\dd\lim_{t\to\infty} s(t)=s_\infty$.
 To discuss the long time behavior of $(u,v)$, we first derive an estimate.

\begin{theo}\lbl{th2.1} \ Let $(u,v,s)$ be the solution of {\rm(NFB)} or {\rm(DFB)}. If
$s_\infty<\infty$, then there exists a constant $K>0$,
such that
 \be
 \|u(t,\cdot),v(t,\cdot)\|_{C^1[0,s(t)]}\leq K,\ \ \ \forall \ t>1.\lbl{2.1}
 \ee
Moveover,
 \be
\lim\limits_{t\to\infty}s'(t)=0.\lbl{2.2}
 \ee
\end{theo}

{\bf Proof}. \ The proof is similar to that of Theorem 4.1 in \cite{WZ} and Lemma 3.3 in \cite{GW}, we omit the details.

\subsection{Vanishing case ($s_\infty<\infty$)}

\begin{theo}\lbl{th2.2} \ Let $(u,v,s)$ be any solution of {\rm(NFB)} or {\rm(DFB)}. If
$s_\infty<\infty$, then
 \be
 \dd\lim_{t\to\infty}\|u(t,\cdot)\|_{C([0,s(t)])}= \dd\lim_{t\to\infty}\|v(t,\cdot)\|_{C([0,s(t)])}=0.\lbl{2.3}
\ee
This result shows that if the species can not spreading into the infinity, they
will die out eventually.
\end{theo}

{\bf Proof}. \  We only prove that $\lim_{t\to\infty}\|u(t,\cdot)\|_{C([0,s(t)])}=0$ for the solution of (NFB), since the rest can be deduced by the similar way. By use of (\ref{1.1}), we have that
 \bes
 u_t\geq u_{xx}+u(1-u-kM), \ \ t>0, \ \ 0<x<s(t).\lbl{2.3a}\ees

The idea used here comes from \cite{WZ}. On the contrary we assume that there exist $\sigma>0$ and $\lk\{(t_j,x_j)\rr\}_{j=1}^{\infty}$, with $0\leq x_j<s(t_j)$ and  $t_j\to\infty$ as $j\to\infty$, such that
 \bes u(t_j,x_j)\geq 3\sigma,\ j=1,2,\cdots.\lbl{2.4}\ees
Since $0\leq x_j<s_\infty$, there are a subsequence of $\{x_j\}$, noted by itself,
and $x_0\in[0,s_\infty]$, such that $x_j\to x_0$ as $j\to\infty$.
We claim that $x_0\in[0,s_\infty)$ . If $x_0=s_\infty$, then $x_j-s(t_j)\to 0$
as $j\to\infty$. By use of the inequality (\ref{2.4}) firstly and the inequality
(\ref{2.1}) secondly, we have that
 \bess
 \lk|\dd\frac{4\sigma}{x_j-s(t_j)}\rr|\leq\lk|\frac{u(t_j,x_j)}{x_j-s(t_j)}\rr|=
 \lk|\frac{u(t_j,x_j)-u(t_j,s(t_j))}{x_j-s(t_j)}\rr|=\lk|u_x(t_j,\bar x_j)\rr|\leq K,\eess
where $\bar x_j\in(x_j,s(t_j))$. It is a contradiction since $x_j-s(t_j)\to 0$.

By use of (\ref{2.1}) and (\ref{2.4}), there exists $\delta>0$ such that $x_0+\delta<s_\infty$ and
 $$u(t_j,x)\geq 2\sigma,\ \ \ \forall\ x\in[x_0, \, x_0+\delta]$$
for all large $j$. Since $s(t_j)\to s_\infty$ as $j\to+\infty$, without loss of generality we may think that $s(t_j)>x_0+\delta$ for all $j$. Let
  \[r_j(t)=x_0+\delta+t-t_j , \ \ \tau_j=\inf\{t>t_j:\, s(t)=r_j(t)\}.\]
Obviously, $x_0+\delta+\tau_j-t_j=r_j(\tau_j)<s_\infty$. Define
 \bess
 y_j(t,x)&=&\dd\frac{2(x-x_0)-(\delta+t-t_j)}{\delta+t-t_j}(\pi-\theta),\\[1mm]
 \underline u_j(t,x)&=&\sigma {\rm e}^{-K(t-t_j)}[\cos y_j(t,x)+\cos
  \theta],\ \ \ (t,x)\in\overline\Omega_{t_j},\eess
where $\theta\ (0<\theta<\pi/8)$ and $K$ are positive constants to be chosen
later, and
   $$\Omega_{t_j}=\{(t,x):t_j<t<\tau_j, \ x_0<x<r_j(t)\}.$$
It is obvious that $\underline u_j(t,x_0)=\underline u_j(t,r_j(t))=0$, and
$|y_j(t,x)|\leq \pi-\theta$ for $(t,x)\in\overline\Omega_{t_j}$, the latter
implies $\underline u_j(t,x)\geq 0$ in $\Omega_{t_j}$.

We want to compare $u(t,x)$ and $\underline u_j(t,x)$ in $\overline\Omega_{t_j}$.
According to (\ref{1.1}), it follows that
 $$u(t,x_0)>0=\underline u_j(t,x_0), \ \ u(t,r_j(t))\geq 0=\underline
u_j(t,r_j(t)), \ \ \forall \  t\in[t_j,\tau_j].$$
On the other hand, it is obvious that
 $$u(t_j,x)\geq 2\sigma>\underline u_j(t_j,x), \ \ \forall \  x\in[x_0,x_0+\delta].$$
Thus, if the positive constants $\theta$ and $K$ can be chosen independent of
$j$ such that
 \bes\underline u_{jt}-\underline u_{jxx}\leq\underline u_j(1-\underline u_j-kM)
  \ \ \ \mbox{in} \ \ \Omega_{t_j},
  \label{2b.8}\ees
it can be deduced that $\underline u_j(t,x)\leq u(t,x)$ for $(t,x)\in\Omega_{t_j}$ by
applying the comparison principle to $w$ and $\underline u_j$ over $\Omega_{t_j}$.
Since $u(\tau_j,s(\tau_j))=0=\underline u_j(\tau_j,r_j(\tau_j))$ and
$s(\tau_j)=r_j(\tau_j)$, it follows that $u_x(\tau_j,s(\tau_j))\leq \underline
u_{jx}(\tau_j,r_j(\tau_j))$. Thanks to $\theta<\pi/8$ and $\delta+\tau_j-t_j<s_\infty$, we have
  $$\underline u_{jx}(\tau_j,r_j(\tau_j))=-\dd\frac{2\sigma(\pi-\theta)}
{\delta+\tau_j-t_j}{\rm e}^{-K(\tau_j-t_j)}\sin(\pi-\theta)
\leq-\dd\frac{7\sigma\pi}{4s_\infty}{\rm e}^{-Ks_\infty}
 \sin \theta.$$
Note the boundary condition $-\mu u_x(\tau_j,s(\tau_j))\leq s'(\tau_j)$, we have
  $$s'(\tau_j)\geq\dd\frac{7\mu\sigma\pi}{4s_\infty}
  {\rm e}^{-Ks_\infty}\sin \theta,$$
which implies $\limsup_{t\to\infty}s'(t)>0$ since $\lim_{j\to\infty}\tau_j\to\infty$. This contradicts to (\ref{2.2}), and so $\lim_{t\to\infty}\|u(t,\cdot)\|_{C([0,s(t)])}=0$.

Now we prove that if $\theta$ and $K$ satisfy
\bes
 &\theta<\dd\frac{\pi}{8}, \ \ \ \sin\theta\leq
\dd\frac{3\pi \delta^2}{4s_\infty^3},&\label{2b.9}\\[1.5mm]
 &K>2\sigma+kM+\left(\dd\frac{2\pi}{\delta}\right)^2+\dd\frac{2\pi s_\infty}{\delta^2(\cos\theta-\cos
2\theta)},&\label{2b.10}
  \ees
then (\ref{2b.8}) holds for all large $j$. And so, $\lim_{t\to\infty}\|u(t,\cdot)\|_{C([0,s(t)])}=0$ is followed. Thanks to $0\leq \underline u_j\leq 2\sigma$ and $\delta+\tau_j-t_j<s_\infty$, a series of computations indicate that, for $(t,x)\in\Omega_{t_j}$,
 \bess
&&\underline u_{jt}-\underline u_{jxx}-\underline u_j(1-\underline u_j-kM)\\[2mm]
&=&-K\underline u_j-\sigma e^{-K(t-t_j)}y_{jt}\sin
y_j+\sigma e^{-K(t-t_j)}y_{jx}^2\cos y_j-\underline
u_j(1-\underline u_j-kM)\\[2mm]
&\leq&\dd\lk(2\sigma +kM+y_{jx}^2-K\rr)\underline
u_j-\sigma e^{-K(t-t_j)}y_{jx}^2\cos
\theta-\sigma e^{-K(t-t_j)}y_{jt}\sin y_j\\[2mm]
&\leq&\dd\left(2\sigma+kM+\left(\frac{2\pi}\delta\right)^2-K\right)\underline u_j
+\sigma{\rm e}^{-K(t-t_j)}\dd\left[\dd\frac{2\pi(x-x_0)}{\delta^2}|\sin y_j|-\left(\dd\frac{2(\pi-\theta)}
{s_\infty}\right)^2\cos\theta\right]\\[1.5mm]
&:=&I.
  \eess
By (\ref{2b.10}), we have $2\sigma+kM+\left(\frac{2\pi}{\delta}\right)^2-K<0$.

Since $-(\pi-\theta)\leq y_j\leq \pi-\theta$ when
$(t,x)\in\overline\Omega_{t_j}$, we can decompose $\Omega_{t_j}=\Omega_{t_j}^1\bigcup
\Omega_{t_j}^2$, where
\bess
  \Omega_{t_j}^1&=&\left\{(t,x)\in\Omega_{t_j}:\, t_j<t<\tau_j,\
 \pi-2\theta<|y_j(t,x)|<\pi-\theta\right\},\\[1mm]
 \Omega_{t_j}^2&=&\left\{(t,x)\in\Omega_{t_j}:\, t_j<t<\tau_j,
 \ |y_j(t,x)|\leq\pi-2\theta\right\}.\eess
It is obvious that $|\sin y_j|\leq \sin 2\theta$ when $(t,x)\in \Omega_{t_j}^1$,
and $\cos y_j\geq -\cos 2\theta$ when $(t,x)\in \Omega_{t_j}^2$. Because of $\theta<\pi/8$, $\underline u_j(t,x)\geq 0$
and $x-x_0\leq s_\infty$ in $\Omega_{t_j}$, in view of (\ref{2b.9}) and (\ref{2b.10}), we conclude
   $$I\leq \sigma {\rm e}^{-K(t-t_j)}\left(\dd\frac{2\pi s_\infty}{\delta^2}\sin 2\theta-\dd\frac{3\pi^2}{s_\infty^2}\cos \theta\right)<0$$
when $(t,x)\in \Omega_{t_j}^1$, and
 $$I\leq \sigma {\rm e}^{-K(t-t_j)}\left(\left[2\sigma+kM
 +\left(\frac{2\pi}{\delta}\right)^2-K\right](\cos\theta-\cos 2\theta)
 +\dd\frac{2\pi s_\infty}{\delta^2}\right)<0$$
when $(t,x)\in \Omega_{t_j}^2$.\ \ \ \fbox{}

\subsection{Spreading case ($s_\infty=\infty$) for the problem (NFB)}

The following theorem is due to Guo and Wu \cite{GW}.

\begin{theo}\lbl{th2.3} Let $(u,v,s)$ be the solution of {\rm(NFB)}. If $s_\infty=\infty$, then for the weakly competition case $0<h,k<1$, we have
\[
 \lim_{t\to\infty}u(t,x)=\frac{1-k}{1-hk},\ \ \ \
\lim_{t\to\infty}v(t,x)=\frac{1-h}{1-hk}\]
uniformly in any compact subset of $[0,\infty)$.
 \end{theo}

In order to investigate the long time behavior of the solution $(u,v)$ to (NFB) for the other cases, we should do some preparation work. We first state two propositions, their proofs are similar to those of Proposition B.1 and Proposition B.2 in \cite{WZ}, respectively.

\begin{prop}\label{p2.1} \ Let $d,\alpha$ and $\beta$ be fixed positive constants. For any given $\vp, L>0$, there exists
$l_\vp>\max\bbb\{L,\frac{\pi}2\sqrt{d/(\alpha\beta)}\bbb\}$ such that, if the continuous and non-negative function $z(t,x)$ satisfies
 $$ \left\{\begin{array}{ll}
  z_t-dz_{xx}\geq\, (\leq)\, \alpha z(\beta-z), \ \ &t>0, \ \ 0<x<l_\vp,\\[2.5mm]
 z_x(t,0)=0, \ \ z(t, l_\vp)\geq\, (=)\, 0, \ \ &t>0,
 \end{array}\right.$$
and $z(0,x)>0$ in $(0,l_\vp)$, then
 \[\liminf_{t\to\infty}z(t,x)>\beta-\vp \ \ \lk(\limsup_{t\to\infty}z(t,x)<\beta+\vp\rr)
 \ \ \ \mbox{uniformly\, on }\, [0,L].\]
  \end{prop}

\begin{prop}\label{p2.2} \ Let $d, \alpha, \beta$ and $K$ be positive constants. For any given $\vp, L>0$, there exists $l_\vp>\max\bbb\{L,\frac{\pi}2\sqrt{d/(\alpha\beta)}\,\bbb\}$ such that, if the continuous and non-negative function $z(t,x)$ satisfies
 $$ \left\{\begin{array}{ll}
  z_t-dz_{xx}\geq\, (\leq)\, \alpha z(\beta-z), \ \ &t>0, \ \ 0<x<l_\vp,\\[2mm]
 z_x(t,0)=0, \ \ z(t, l_\vp)\geq\, (\leq)\, K,\ \ &t>0,
 \end{array}\right.$$
and $z(0,x)>0$ in $(0,l_\vp)$, then we have
  \[\liminf_{t\to\infty}z(t,x)\geq \beta-\vp \ \ \lk(\limsup_{t\to\infty}z(t,x)<\beta+\vp\rr)
  \ \ \ \mbox{uniformly\, on }\, [0,L].\]
 \end{prop}

Now we give another proposition.

\begin{prop}\label{p2.3} \ Let $d, \lambda, K>0$ and $\zeta\geq 0$ be fixed constants. For any given $\vp, L>0$, there exists $l_\vp\gg 1$ such that, if the continuous and non-negative function $z(t,x)$ satisfies
 $$ \left\{\begin{array}{ll}
  z_t-dz_{xx}\leq -\la z(\zeta+z), \ \ &t>0, \ \ 0<x<l_\vp,\\[2mm]
 z_x(t,0)=0, \ \ z(t, l_\vp)\leq K,\ \ &t>0,\\[2mm]
 0<z(0,x)\leq K, & x\in (0,l_\vp),
 \end{array}\right.$$
then we have
 \[\limsup_{t\to\infty}z(t,x)<\vp \ \ \ \mbox{uniformly\, on }\, [0,L].\]
 \end{prop}

{\bf Proof}. Let $\theta\in(0,K)$ and $y_\theta(x)$ be the unique positive solution of
 \bedd
 dy''=\la y(\zeta+y), \ \ x>0,\\
 y(0)=\theta, \ \ y'(0)=0.
 \nm\eedd
Then $y_\theta(x)$ is increasing in $x$, and there exist two positive constants $l_\theta$ and $L_\theta$, with $l_\theta<L_\theta$ and $\lim_{\theta\to 0}l_\theta=\lim_{\theta\to 0}L_\theta=\infty$, such that $y_\theta(l_\theta)=K$ and $\lim_{x\to L_\theta}y_\theta(x)=\infty$. Moreover, $\lim_{\theta\to 0}y_\theta(x)=0$ uniformly in any compact subset of $[0,\infty)$. For any given $\ep>0$, there is a $\theta>0$ such that
 \bes
 l_\theta>L, \ \ y_\theta(x)<\ep\ff x\in[0,L].\lbl{2.9}\ees

Let $l_\ep\in(l_\theta,L_\theta)$ be fixed, and $w(t,x)$ be the unique positive solution of
 \bedd
 w_t-dw_{xx}=-\la w(\zeta+w), \ &t>0, \ \ 0<x<l_\ep,\\\medskip
 w_x(t,0)=0, \ \ w(t,l_\ep)= K,\ \ &t>0, \\
 w(0,x)=K, &0\leq x\leq l_\ep.
 \nm\eedd
It is easy to know that the limit $\lim_{t\to\infty}w(t,x)=W(x)$ exists and is a positive solution of the following boundary value problem
 \bedd
 dW''=\la W(\zeta+W), \ \ \ 0<x<l_\ep,\\
 W'(0)=0, \ \ W(l_\ep)= K.
 \nonumber\eedd
Obviously, $W(x)$ is increasing in $x$. Since $y_\theta(l_\ep)>K$, we get $W(0)<y_\theta(0)$ by the comparison principle for the initial value problem of ODE. Hence, applying such comparison principle once again, one has $W(x)\leq y_\theta(x)$ in $[0,l_\ep]$, and so in $[0,L]$. Consequently, by (\ref{2.9}), $W(x)<\ep$ in $[0,L]$. It is deduced that
 \[\limsup_{t\to\infty}w(t,x)<\vp \ \ \ \mbox{uniformly\, on }\, [0,L].\]
Applying the comparison principle for parabolic equations to $z(t,x)$ and $w(t,x)$, we can get the desired conclusion.  \ \ \ \fbox{}

\begin{theo}\lbl{tb2.4}\, Let $(u,v,s)$ be the solution of {\rm(NFB)}, and assume that $s_\infty=\infty$. We have the following conclusions:

{\rm(i)}\, If $0<h<1\leq k$, then
 \bes \lim_{t\to\infty}u(t,x)=0,\ \ \lim_{t\to\infty}v(t,x)=1\ \ {\rm uniformly\ in\ any\ compact\ subset\ of}\ [0,\infty).\lbl{2.10}\ees

{\rm(ii)}\, If $0<k<1\leq h$, then
 \bes\lim_{t\to\infty}u(t,x)=1,\ \ \lim_{t\to\infty}v(t,x)=0\ {\rm uniformly\ in\ any\ compact\ subset\ of}\ [0,\infty).\lbl{2.11}\ees
 \end{theo}

\begin{remark}\lbl{r2.1} By the conclusions of Theorems $\ref{th2.3}$ and $\ref{tb2.4}$ we see that, if $s_\infty=\infty$, the dynamic behaviors of {\rm(NFB)} are same as those of the positive solutions of ODE system
 \bedd
 u'=u(1-u-kv), \ \ &t>0,\\
 v'=rv(1-v-hu),&t>0
 \lbl{2.12}\eedd
for the following cases: if $0<h,k<1$, then $\lim_{t\to\infty}(u,v)(t)=\lk(\frac{1-k}{1-hk}, \frac{1-h}{1-hk}\rr)$; if $0<k<1\leq h$, then $\lim_{t\to\infty}(u,v)(t)=(1,0)$; if $0<h<1\leq k$, then $\lim_{t\to\infty}(u,v)(t)=(0,1)$.
 \end{remark}

The proof of Theorem \ref{tb2.4} will be divided into three lemmas.

\begin{lem}\lbl{l2.1} Assume that $0<h<1\leq k$. If $k(1-h)\geq 1$, then $(\ref{2.10})$ holds.
\end{lem}

{\bf Proof}. By our assumption, $k>1$. Let $M$ be as in (\ref{1.1}). For any given $L>0$ and $0<\vp\ll 1$, let $l_\vp$ be given by Proposition \ref{p2.2} with $d=1$, $\beta=\alpha=1$ and $K=M$. In view of $s_\infty=\infty$, there exists $T_1>0$ such that
 \[ s(t)>l_\vp, \ \ \ \forall \ t\geq T_1, \ \ x\in[0,l_\vp].\]
Notice that $v> 0$ in $[0,l_\ep]$, we see that $u$ satisfies
 \bess
 \left\{\begin{array}{lll}
  u_t-u_{xx}\leq u(1-u),\ &t\geq T_1, \ \ x\in[0,l_\vp],\\[2mm]
  u_x(t,0)=0, \ \  u(t,l_\vp)\leq M,\ \ &t\geq T_1.
 \end{array}\right.
 \eess
Since $u(T_1,x)>0$ in $[0,l_\vp]$, applying Proposition \ref{p2.2}, it arrives at
 \[\limsup_{t\to\infty}u(t,x)\leq 1+\vp\ \ \ \mbox{uniformly\, on }\, [0,L].\]
By the arbitrariness of $\vp$ and $L$,
  \bes\limsup_{t\to\infty}u(t,x)\leq 1:=\bar{u}_1\ \ \mbox{uniformly\, on\, the\, compact\,  subset\, of } \, [0,\infty). \label{2.13}\ees

For any given $L>0$, $0<\delta\ll
1$ with $h(1+\delta)<1$ and $0<\vp\ll 1$, let $l_\vp$ be given by Proposition \ref{p2.2} with $d=D$, $\alpha=r$ and $\beta=1-h\lk(\bar u_1+\delta\rr)$. In view of (\ref{2.13}), there exists $T_2>0$ such that $u(t,x)\leq \bar u_1+\delta$ for all
$t\geq T_2$ and $x\in[0,l_\vp]$. Therefore, $v$ satisfies
 \bess
 \left\{\begin{array}{lll}
  v_t-Dv_{xx}\geq rv\lk[1-v-h\lk(\bar u_1+\delta\rr)\rr],
  \ &t\geq T_2, \ \ x\in[0,l_\vp],\\[2mm]
  v_x(t,0)=0, \ \ v(t, l_\vp)\geq 0,\ \ &t\geq T_2.
 \end{array}\right.
 \eess
As $v(T_2,x)>0$ in $[0,l_\vp]$, in view of Proposition \ref{p2.1} it yields
 \[\liminf_{t\to\infty}v(t,x)\geq 1-h\big(\bar u_1+\delta\big)-\vp\ \ \ \mbox{uniformly\, on }\, [0,L].\]
By the arbitrariness of $\vp$, $L$ and $\delta$ we have
\bes\liminf_{t\to\infty}v(t,x)\geq 1-h\bar u_1:=\ud{v}_1\ \ \mbox{uniformly\, on\, the\, compact\,  subset\, of } \, [0,\infty). \label{2.14}\ees

When $k(1-h)>1$, we choose $0<\delta\ll 1$ such that $1-k(\ud{v}_1-\delta)<0$.
For any given $L>0$ and $0<\vp\ll 1$, let $l_\vp$ be given by
Proposition \ref{p2.3} with $d=1$, $\zeta=k(\ud v_1-\delta)-1>0$, $\la=1$ and
$K=M$. Taking into account (\ref{2.14}) and $s_\infty=\infty$, there is $T_3>0$ such that
 \[v(t,x)\geq \ud v_1-\delta,\ \ s(t)>l_\vp, \ \ \
  \forall \ t\geq T_3, \ \ x\in[0,l_\vp].\]
It yields that $u$ satisfies
 \bess
 \left\{\begin{array}{lll}
  u_t-u_{xx}\leq u\lk[1-u-k\lk(\ud{v}_1-\delta\rr)\rr],\ &t\geq T_3,
   \ \ x\in[0,l_\vp],\\[2mm]
  u_x(t,0)=0,\ \ \ u(t,l_\vp)\leq K,\ \ &t\geq T_3.
 \end{array}\right.
 \eess
By use of Proposition \ref{p2.3} we have
 \[\limsup_{t\to\infty}u(t,x)<\vp \ \ \ \mbox{uniformly\, on }\, [0,L].\]
Note that $u\geq 0$, by the arbitrariness of $\vp$ and $L$, we have
\bes\lim_{t\to\infty}u(t,x)=0\ \ \mbox{uniformly\, on\, the\, compact\,  subset\, of } \, [0,\infty). \label{2.15}\ees

When $k(1-h)=1$. For any given $L>0$, $0<\delta\ll 1$ and $0<\vp\ll 1$, let $l_\vp$ be given by
Proposition \ref{p2.2} with $d=\alpha=1$, $\beta=k\delta$ and
$K=M$. Taking into account (\ref{2.14}) and $s_\infty=\infty$, there is $T_4>0$ such that
 \[v(t,x)\geq \ud v_1-\delta,\ \ s(t)>l_\vp, \ \ \
  \forall \ t\geq T_4, \ \ x\in[0,l_\vp].\]
It yields that $u$ satisfies
 \bess
 \left\{\begin{array}{lll}
  u_t-u_{xx}\leq u\lk[1-u-k\lk(\ud{v}_1-\delta\rr)\rr]=u(k\delta-u),\ &t\geq T_4,
   \ \ x\in[0,l_\vp],\\[2mm]
  u_x(t,0)=0,\ \ \ u(t,l_\vp)\leq K,\ \ &t\geq T_4.
 \end{array}\right.
 \eess
By use of Proposition \ref{p2.2} we have
 \[\limsup_{t\to\infty}u(t,x)<k\delta+\vp \ \ \ \mbox{uniformly\, on }\, [0,L].\]
Note that $u\geq 0$, by the arbitrariness of $\vp$, $\delta$ and $L$, we have that (\ref{2.15}) is still true.

For any given $L>0$, $0<\delta\ll 1$ and $0<\vp\ll 1$, let $l_\vp$ be given by
Proposition \ref{p2.1} with $d=D$, $\alpha=r$ and $\beta=1-h\delta$. According to (\ref{2.15}), there is $T_5>0$ such that $u(t,x)\leq\delta$ for
all $t\geq T_5$ and $x\in[0,l_\vp]$. Hence, $v$ satisfies
 \bess
 \left\{\begin{array}{lll}
  v_t-Dv_{xx}\geq rv\lk[1-v-h\delta\rr],\ &t\geq T_5,
   \ \ x\in[0,l_\vp],\\[2mm]
 v_x(t,0)=0, \ \  v(t,l_\vp)\geq 0,\ \ &t\geq T_5.
 \end{array}\right.
 \eess
Similar to the above,
 \bes\liminf_{t\to\infty}v(t,x)\geq 1\ \ \mbox{uniformly\, on\, the\,
compact\,  subset\, of } \,[0,\infty). \lbl{2.16}\ees

Let $w(t,x)$ be the unique positive solution of
 \bess
 \left\{\begin{array}{lll}
  w_t-Dw_{xx}= rw(1-w),\ &t>0, \ \ 0<x<\infty,\\[1mm]
 w_x(t,0)=0, \ \ &t>0,\\[1mm]
  w(0,x)=\phi(x) & x\geq 0,
  \end{array}\right.
 \eess
where,
 \[\phi(x)=\left\{\begin{array}{ll}
  u_0(x), \ \ & 0\leq x\leq s_0,\\[2mm]
  0, \ \ & x\geq s_0.
  \end{array}\right.\]
It is well known that $\lim_{t\to\infty}w(t,x)=1$ uniformly on the compact subset of $[0,\infty)$ since $u_0(x)>0$ in $(0,s_0)$. Furthermore, the comparison principle yields $v\leq w$. Consequently,
$\limsup_{t\to\infty}v(t,x)\leq 1$ uniformly on the compact subset of $[0,\infty)$.
Remember (\ref{2.16}), one has $\lim_{t\to\infty}v(t,x) =1$ uniformly on the compact subset of $[0,\infty)$. Thanks to the limit (\ref{2.15}), we obtain (\ref{2.10}). The proof is complete. \ \ \ \fbox{}

\begin{lem}\lbl{l2.2}  Assume that $0<h<1\leq k$. If $k(1-h)<1$, then $(\ref{2.10})$ holds.
 \end{lem}

{\bf Proof}.\, Follow the proof of Lemma \ref{l2.1}, by our assumption, $1-k\ud{v}_1>0$, where $\ud v_1$ is given in (\ref{2.14}). For any given $L>0$, $0<\delta\ll 1$ and $0<\vp\ll 1$, let $l_\vp$ be given by Proposition \ref{p2.2} with $d=1$, $\beta=1-k\lk(\ud v_1-\delta\rr)$, $\alpha=1$ and $K=M$. Taking into account (\ref{2.14}) and $s_\infty=\infty$, there is $T_6>0$ such that
 \[v(t,x)\geq \ud v_1-\delta,\ \ s(t)>l_\vp, \ \ \
  \forall \ t\geq T_6, \ \ x\in[0,l_\vp].\]
It yields that $u$ satisfies
 \bess
 \left\{\begin{array}{lll}
  u_t-u_{xx}\leq u\lk[1-u-k\lk(\ud{v}_1-\delta\rr)\rr],\ &t\geq T_6,
   \ \ x\in[0,l_\vp],\\[2mm]
  u_x(t,0)=0,\ \ \ u(t,l_\vp)\leq K,\ \ &t\geq T_6.
 \end{array}\right.
 \eess
Similar to the above, we can get
 \bes\limsup_{t\to\infty}u(t,x)\leq 1-k\ud v_1:=\bar u_2\ \
 \mbox{uniformly\, on\, the\, compact\,  subset\, of } \, [0,\infty). \label{2.17}\ees
By our assumption, $\bar u_2=1-k(1-h)>0$.

Since $h<1$, we see that
 \[\underline{v}_2:=1-h\bar u_2=(1-h)(1+hk)>0. \]
For any given $L>0$, $0<\delta\ll 1$ and $0<\vp\ll 1$, let $l_\vp$ be given by
Proposition \ref{p2.1} with $d=D$, $\alpha=r$ and $\beta=1-h\lk(\bar u_2+\delta\rr)$.
According to (\ref{2.17}), there is $T_7>0$ such that $u(t,x)\leq\bar u_2+\delta$ for
all $t\geq T_7$ and $x\in[0,l_\vp]$. Hence, $v$ satisfies
 \bess
 \left\{\begin{array}{lll}
  v_t-Dv_{xx}\geq rv\lk[1-v-h\lk(\bar u_2+\delta\rr)\rr],\ &t\geq T_7,
   \ \ x\in[0,l_\vp],\\[2mm]
 v_x(t,0)=0, \ \  v(t,l_\vp)\geq 0,\ \ &t\geq T_7.
 \end{array}\right.
 \eess
Similar to the above,
 \bess\liminf_{t\to\infty}v(t,x)\geq \underline{v}_2\ \ \mbox{uniformly\, on\, the\,
compact\,  subset\, of } \, [0,\infty). \eess

If $k\ud{v}_2\geq 1$, similar to the proof of Lemma \ref{l2.1}, we can get (\ref{2.10}).

If $k\ud{v}_2<1$, we continue the above iterative process and get
\bess\liminf_{t\to\infty}u(t,x)\leq
1-k\underline{v}_2:=\bar{u}_3\ \ \mbox{uniformly\, on\, the\,
compact\,  subset\, of } \, [0,\infty),\\
 \liminf_{t\to\infty}v(t,x)\geq 1-h\bar u_3:=\ud{v}_3\ \ \mbox{uniformly\, on\, the\,
compact\,  subset\, of } \, [0,\infty).
\eess

For $j\geq 3$, if $k\ud v_i<1$ for all $2\leq i\leq j$, we define
 \[\bar{u}_{j+1}=1-k\ud v_j, \ \ \ \ud{v}_{j+1}=1-h\bar u_{j+1}.\]
Then $\bar u_j,\ud v_j>0$ for all such $j$ since $h<1$.

If there is a first $j\geq 3$ such that $k\ud v_j\geq 1$,  similar to the proof of Lemma \ref{l2.1}, we can get (\ref{2.10}).

Now we assume that $k\ud v_j<1$ for all $j\geq 3$. Repeating the above procedure, we have that, for all $j\geq 3$,
\bedd\liminf_{t\to\infty}u(t,x)\leq\bar{u}_j\ \ \mbox{uniformly\, on\, the\,
compact\,  subset\, of } \, [0,\infty),\\
 \liminf_{t\to\infty}v(t,x)\geq \ud{v}_j\ \ \mbox{uniformly\, on\, the\,
compact\,  subset\, of } \, [0,\infty).
  \lbl{2.18}\eedd
Carefully calculation gives
 \[\ud v_j=(1-h)\bbb(1+\sigma+\sigma^2+\sigma^3+\cdots+\sigma^{j-1}\bbb), \ \ \ \ j\geq 3,\]
where $\sigma=hk$. By our assumption, $k\ud v_j<1$ for all $j\geq 3$, which implies $hk=\sigma<1$. Therefore,
 \[\lim_{j\to\infty}\ud v_j=\frac{1-h}{1-hk}, \ \ \mbox{and consequently} \  \lim_{j\to\infty}\bar u_j=\frac{1-k}{1-hk}.\]
Since $\bar u_j,\ud v_j>0$, we conclude that $k\leq 1$. By the assumption of this lemma, it yield $k=1$. Thus, we have $\lim_{j\to\infty}\bar u_j=0$. Remember the first limit of (\ref{2.18}) and the fact that $u(t,x)\geq 0$, we obtain (\ref{2.15}). Similar to the last part of the proof of Lemma \ref{l2.1}, it follows that (\ref{2.10}) holds. The proof is complete. \ \ \ \fbox{}

Similar to the above we can prove

\begin{lem}\lbl{l2.3} If $0<k<1\leq h$, then $(\ref{2.11})$ holds.
\end{lem}

The conclusions of Theorem $\ref{tb2.4}$ can be followed by Lemmas \ref{l2.1}--\ref{l2.3}.

\subsection{Spreading case ($s_\infty=\infty$) for the problem (DFB)}

We first state two propositions.

 \begin{prop}\lbl{p2.4}{\rm(\cite{W})} \ Let $d$ and $\lambda$ be positive constants. Assume that $f$ satisfies
 \bess
 f\in C^\alpha_{\rm loc}([0,\infty)) \ \ \mbox{with} \ 0<\alpha<1, \ \
 \inf_{x\geq 0}f(x):=f_0>0, \ \ \|f\|_\infty<\infty.
 \eess
Then the problem
 \bess\left\{\begin{array}{ll}\medskip
 -du''=u\big(f(x)-\lambda u\big), \ \ 0<x<\infty,\\
 u(0)=0,
 \end{array}\right.\eess
has a unique positive solution $u(x)$. Furthermore,

{\rm(i)}\, if $f(x)$ is increasing in $x$, so is $u(x)$ and $\lim\limits_{x\to\infty}u(x)=
\frac 1\lambda\lim\limits_{x\to\infty}f(x)$;

\vskip 4pt {\rm(ii)}\, if $f(x)$ is decreasing in $x$, then either $u(x)$ is increasing in $x$, or there exists $x_0>0$ such that $u(x)$ is increasing in $(0,x_0)$ and $u(x)$ is decreasing in $(x_0,\infty)$. Therefore, $\lim\limits_{x\to\infty}u(x)=\frac 1\lambda\lim\limits_{x\to\infty}f(x)$.
\end{prop}

Let $d$ and $\lambda$ be as above. Assume that $f_i\in C^\alpha_{\rm loc}([0,\infty))$ and satisfies $0<\inf_{x\geq 0}f_i(x)\leq\|f_i\|_\infty<\infty$ for $i=1,2$. By Proposition  \ref{p2.4}, the problem
 \bes\left\{\begin{array}{ll}\medskip
 -du''=u\big(f_i(x)-\lambda u\big), \ \ 0<x<\infty,\\
 u(0)=0
 \nm\end{array}\right.\ees
has a positive solution, denoted by $u_i$.

\begin{prop}{\rm(}Comparison principle, {\rm\cite{W})}\lbl{p2.5} \ Under the above conditions,
if $f_1(x)\leq f_2(x)$ for all $x\geq 0$, then we have that
 \[u_1(x)\leq u_2(x), \ \ \ \forall \ x\geq 0.\]
 \end{prop}

 \begin{theo}\lbl{th2.5} \ Assume that $0<h,k<1$. Then the problem
 \bedd
 -u''=u(1-u-kv), &0<x<\infty,\\\medskip
 -Dv''=rv(1-v-hu),\ \ &0<x<\infty, \\
 u(0)=v(0)=0
 \lbl{2.19}\eedd
has a positive solution. Moreover, any positive solution $(u,v)$ of {\rm(\ref{2.19})} satisfies
 \bes
 \ud u(x)\leq u(x)\leq\bar u(x), \ \ \ \ud v(x)\leq v(x)\leq\bar v(x), \ \ \
 \forall \ x\geq 0,\lbl{2.20}\ees
where $\bar u,\,\bar v,\,\ud u$ and $\ud v$ will be given in the proof.
\end{theo}

{\bf Proof}. {\it Step 1}: The construction of $\ud u$, $\ud v$, $\bar u$ and $\bar v$.

Let $\bar u$ be the unique positive solution of
\[
\left\{\begin{array}{l}
 -u''=u(1-u), \ \ 0<x<\infty,\\[2mm]
 u(0)=0.
\end{array}\right. \]
Then
 \[\bar u'(x)>0, \ \ \mbox{and} \ \ \lim_{x\to\infty}\bar u(x)=1.\]
Since $h<1$, in view of Proposition \ref{p2.4}, the problem
 \bedd
 -Dv''=rv\big[1-v-h\bar u(x)\big],\ \ &0<x<\infty, \\
 v(0)=0
 \nm\eedd
has a unique positive solution, denoted by $\ud v(x)$, and $\ud v(x)\leq 1-h$. Let $\bar v$ be the unique positive solution of
 \bedd
 -Dv''=rv\big(1-v(x)\big), \ \ 0<x<\infty,\\
 v(0)=0.
 \nm\eedd
Then
 \[\bar v'(x)>0, \ \ \mbox{and} \ \ \lim_{x\to\infty}\bar v(x)=1.\]
Thanks $k<1$, by virtue of Proposition \ref{p2.4}, the problem
 \bedd
 -u''=u\big[1-u-k\bar v(x)\big],\ \ &0<x<\infty, \\
 u(0)=0
 \nm\eedd
has a unique positive solution, denoted by $\ud u(x)$, and $\ud u(x)\leq 1-k$.

Applying Proposition \ref{p2.5}, we have that $\ud u(x)\leq \bar u(x)$ and $\ud v(x)\leq \bar v(x)$ for all $x\geq 0$.

{\it Step 2}: Existence of positive solutions.

The conclusion of Step $1$ show that $\ud u$, $\ud v$, $\bar u$ and $\bar v$ are the coupled ordered lower and upper solutions of (\ref{2.19}). For any given $l>0$, it is obvious that $\ud u$, $\ud v$, $\bar u$ and $\bar v$ are also the coupled ordered lower and upper solutions of the following problem
 \bedd
 -u''=u(1-u-kv), &0<x<l,\\\medskip
 -Dv''=rv(1-v-hu),\ \ &0<x<l, \\
 u(0)=v(0)=0, \ \ u(l)=\bar u(l), \ \ v(l)=\bar v(l).
 \lbl{2.21}\eedd
By the standard upper and lower solutions method we have that the problem (\ref{2.21}) has a least one positive solution, denoted by $(u_l,v_l)$ and
 \bes\ud u(x)\leq u_l(x)\leq\bar u(x), \ \ \ \ud v(x)\leq v_l(x)\leq\bar v(x), \ \ \
 \forall \ 0\leq x\leq l.\nm\ees
Applying the local estimation and compactness argument, it can be concluded that there exists a pair $(u,v)$, such that $(u_l,v_l)\longrightarrow (u,v)$ in $\lk[C^2_{\rm loc}([0,\infty))\rr]^2$, and $(u, v)$ solves (\ref{2.19}).

By use of Proposition \ref{p2.5}, it is easy to see that if $0<k,h<1$ then any positive solution $(u,v)$ of (\ref{2.19}) satisfies (\ref{2.20}).\ \ \ \fbox{}

Similar to the proof of Theorem 3.5 in \cite{W}, we can prove the following

\begin{theo}\lbl{th2.6} \ Assume that $s(\infty)=\infty$. If $0<k,h<1$, then the solution $(u(t,x),v(t,x))$ of {\rm(DFB)} satisfies
 \bess
 &\dd\liminf_{t\to\infty}u(t,x)\geq\ud u(x), \ \ \limsup_{t\to\infty}u(t,x)\leq\bar u(x)\ \ \ \mbox{uniformly in any compact subset of } \, [0, \infty),\qquad&\\
 &\dd\liminf_{t\to\infty}v(t,x)\geq\ud v(x), \ \ \limsup_{t\to\infty}v(t,x)\leq\bar v(x)\ \ \ \mbox{uniformly in any compact subset of } \, [0, \infty),\qquad&
 \eess
where $\bar u,\,\bar v,\,\ud u$ and $\ud v$ are given in the proof of Theorem $\ref{th2.5}$.
\end{theo}

\section{The criteria governing spreading and vanishing}
\setcounter{equation}{0}

In this section we first state a comparison principle which can be used to determine the criteria governing spreading and vanishing.

\begin{lem} $($Comparison principle$)$\label{le3.1} \
Let $\bar s\in C^1([0,T])$ and $\bar s(t)>0$ in $[0,\infty)$. Let $\bar u, \bar v\in
C(\overline{O})\bigcap C^{1,2}(O)$ with $O=\{(t,x): t>0,\, 0<x<\bar s(t)\}$. Assume that $(\bar u,\bar v,
\bar s)$ satisfies
\be
 \left\{\begin{array}{ll}
  \bar u_t-\bar u_{xx}\geq\bar u(1-\bar u),\ \ &t>0,\ \ 0<x<\bar s(t),\\[2mm]
 \bar v_t-D\bar v_{xx}\geq r\bar v(1-\bar v),&t>0,\ \ 0<x<\bar s(t),\\[2mm]
 \bar u_x(t,0)\leq 0,\ \ \bar v_x(t,0)\leq 0,&t>0,\\[2mm]
 \bar u(t,\bar s(t))=\bar v(t,\bar s(t))=0,&t>0,\\[2mm]
 \bar s'(t)\geq-\mu[\bar u_x(t,\bar s(t))+\rho \bar v_x(t,\bar s(t))],\ \ &t>0.
 \end{array}\right.\lbl{3.1}
\ee
If $\bar s(0)\geq s_0,\ \bar u(0,x)\geq 0,\ \bar v(0,x)\geq 0$ on $[0,\bar s(0)]$, and
 $u_0(x)\leq\bar u(0,x),\ v_0(x)\leq\bar v(0,x)$ on $[0,s_0]$,
then the solution $(u,v,s)$ of {\rm (NFB)} satisfies
\[
  s(t)\leq\bar s(t) \ {\rm on}\ [0,\infty),\ \ u(t,x)\leq\bar u(t,x),\ \  v(t,x)\leq\bar v(t,x)\ \ {\rm on}\ \ \overline\Omega,\]
where $\Omega=\{(t,x):\ t>0,\ 0<x<s(t)\}$.

If, in $(\ref{3.1})$, the conditions $\bar u_x(t,0)\leq 0$ and $\bar v_x(t,0)\leq 0$ are replaced by
$\bar u(t,0)\geq 0$ and $\bar v(t,0)\geq 0$, then the conclusion still holds for the solution of {\rm(DFB)}.
\end{lem}

{\bf Proof}. \ The proof is same as that of \cite[Lemma 4.1]{W} (see also the argument of \cite[Lemma 5.1]{GW}), we omit the details. \ \ \ $\Box$

\vskip 4pt The pair $(\bar u,\bar v,\bar s)$ in Lemma \ref{le3.1} is usually called an
upper solution of (NFB) or (DFB).

We next give a necessary condition of vanishing.

\begin{theo}\lbl{th3.1} If $s_\infty<\infty$, then $s_\infty\leq \frac\pi 2\min\lk\{1,\sqrt{ D/r}\rr\}$ for the problem {\rm (NFB)},
and $s_\infty\leq \pi\min\lk\{1,\sqrt{ D/r}\rr\}$ for the problem {\rm (DFB)}
\end{theo}

{\bf Proof}. The proof is similar to that of \cite[Lemma 5.1]{W}, we omit the details.
\ \ \ $\Box$

\vskip 2pt Define
 \[\Lambda=\left\{\begin{array}{l}
\frac\pi 2\min\lk\{1,\sqrt{D/r}\rr\}\ {\rm for\ the\ problem\ (NFB)},\\[3mm]
\pi\min\lk\{1,\sqrt{D/r}\rr\}\ {\rm for\ the\ problem\ (DFB)}.
  \end{array}\right.\]
Then by Theorem \ref{th3.1} and (\ref{1.1}), $s_0\geq \Lambda$
implies $s_\infty=\infty$.

In the following, with the parameters $s_0$ satisfying $s_0<\Lambda$ and $(u_0,v_0)$ fixed, let us discuss the effect of the coefficient $\mu$ on the spreading and vanishing. As a first step, when $\mu$ is sufficiently large, we have

\begin{lem}\lbl{le3.3} \ Suppose that $s_0<\Lambda$. Then for both problems {\rm (NFB)} and {\rm (DFB)}, there exists $\mu^0>0$ depending on $(u_0,v_0,s_0)$ such that $s_\infty=\infty$ if $\mu\geq\mu^0$.
\end{lem}

{\bf Proof}. \ The idea of this proof comes from \cite[Lemma 3.6]{PZ}. We only deal with the problem (NFB), since the problem (DFB) can be treated by the similar way.

We see from (\ref{1.1}) that there exists a constant $\delta^*>0$ satisfying
  \[u(1-u-kv)\geq-\delta^*u,\ \ rv(1-v-hu)\geq-\delta^*v\]
for all $u,v\in[0,M]$. We next consider the auxiliary free boundary problem
\begin{equation}\left\{\begin{array}{lll}
 w_t-w_{xx}=-\delta^*w, &t>0, \ \ 0<x<r(t),\\[1mm]
 z_t-Dz_{xx}=-\delta^*z,&t>0, \ \ 0<x<r(t),\\[1mm]
 w_x=z_x=0, \ &t>0, \ \ x=0,\\[1mm]
w=z=0, \ \ r'(t)=-\mu(w_x+\rho z_x),\ \ &t>0, \ \ x=r(t),\\[1mm]
 w(0,x)=u_0(x), \ \ z(0,x)=v_0(x),& 0\leq x\leq s_0,\\[1mm]
  r(0)=s_0.
 \end{array}\right.\lbl{3.2}
 \end{equation}
Arguing as in proving the existence and uniqueness of the solution to (NFB),
one will easily see that (\ref{3.2}) also admits a unique solution $(w,z,r)$ which is well
defined for all $t>0$. Moreover, due to the Hopf boundary lemma, $r'(t)>0$ for
$t>0$. To stress the dependence of the solutions on the parameter $\mu$, in the sequel, we always write $(u^\mu,v^\mu,s^\mu)$ and $(w^\mu,z^\mu,r^\mu)$ instead of $(u,v,s)$ and $(w,z,r)$. Similar to the proof of \cite[Lemma 4.1]{W}, we have
\be
u^\mu(t,x)\geq w^\mu(t,x),\ \ v^\mu(t,x)\geq z^\mu(t,x),\ \ s^\mu(t)\geq r^\mu(t),\ \forall \  t\geq 0,\ x\in[0,r^\mu(t)].\lbl{3.3}
\ee
In what follows, we are going to prove that for all large $\mu$,
\be
r^\mu(2)\geq 2\Lambda.\lbl{3.4}
\ee

To the end, we first choose a smooth function $\underline r(t)$ with $\underline r(0)=s_0/2$, $\underline r'(t)>0$ and $\underline r(2)=2\Lambda$. We then consider the following initial-boundary value problem
\begin{equation}\left\{\begin{array}{lll}
 \underline w_t-\underline w_{xx}=-\delta^*\underline w, &t>0,\ \ 0<x<\underline r(t),\\[1mm]
 \underline z_t-D\underline z_{xx}=-\delta^*\underline z,&t>0,\ \ 0<x<\underline r(t),\\[1mm]
 \underline w_x(t,0)=\underline z_x(t,0)=0,\ &t>0,\\[1mm]
 \underline w(t,\underline r(t))=\underline z(t,\underline r(t))=0,&t>0,\\[1mm]
 \underline w(0,x)=\underline w_0(x), \ \ \underline z(0,x)=\underline z_0(x),\ \ & 0\leq x\leq s_0/2.
 \end{array}\right.\lbl{3.5}
 \end{equation}
Here, for the smooth initial value $(\underline w_0,\underline z_0)$, we require
 \bedd
0<\underline w_0(x)\leq u_0(x) \ \ {\rm on} \ [0,s_0/2],\ \ \underline w_0'(0)=\underline w_0(s_0/2)=0,\ \ \underline w_0'(s_0/2)<0,\\
0<\underline z_0(x)\leq v_0(x) \ \ {\rm on} \ [0,s_0/2],\ \ \underline z_0'(0)=\underline z_0(s_0/2)=0,\ \ \underline z_0'(s_0/2)<0.
 \lbl{3.6}
 \eedd
The standard theory for parabolic equations ensures that (\ref{3.5}) has a unique positive solution $(\underline w,\underline z)$, and $\underline w_x(t,\underline r(t))<0$, $\underline z_x(t,\underline r(t))<0$ for all $t\in[0,2]$ due to the Hopf boundary lemma. According to our choice of $\underline r(t)$ and $(\underline w_0(x),\underline z_0(x))$, there is a constant $\mu^0>0$ such that, for all $\mu\geq \mu^0$,
\be
\underline r'(t)\leq -\mu\big[\underline w_x(t,\underline r(t))+\rho\underline z_x(t,\underline r(t))\big]\ff 0\leq t\leq 2.\lbl{3.7}
\ee
On the other hand, for system (\ref{3.2}), we can establish the comparison principle
analogous with lower solution to Lemma \ref{le3.1} by the same argument. Thus, note that $\underline r(0)=s_0/2<r^\mu(0)$, it follows from (\ref{3.2}), (\ref{3.5}), (\ref{3.6}) and (\ref{3.7}) that
 \[w^\mu(t,x)\geq \underline w(t,x),\ \ z^\mu(t,x)\geq \underline z(t,x),\ \ r^\mu(t)\geq \underline r(t),\ \ \forall\ t\in[0,2],\ x\in[0,\,\underline r(t)].\]
Which particularly implies $r^\mu(2)\geq \underline r(2)=2\Lambda$, and so (\ref{3.4}) holds true. Hence, in view of (\ref{3.3}) and (\ref{3.4}), we find
\[s_\infty=\lim\limits_{t\to\infty}s^\mu(t)>s^\mu(2)\geq2\Lambda.\]
This, together with Theorem \ref{th3.1}, yields the desired result. \ \ \fbox{}

Secondly, when $s_0<\Lambda$, Guo and Wu \cite{GW} have proved that $s_\infty<\infty$ if $\mu$ is small enough for the problem (NFB).
We give the following assertion.

\begin{lem}\lbl{le3.4} \ Assume that $s_0<\Lambda$. Then
there exists $\mu_0>0$ for the problems {\rm (DFB)}, depending also on $(u_0,v_0,s_0)$, such that
$s_\infty<\infty$ when $\mu\leq\mu_0$.
\end{lem}

{\bf Proof}. \
 We shall use the argument from Ricci and Tarzia \cite{RT} to construct the suitable upper
solutions and use Lemma 3.1 to derive the desired conclusion. We adopt the following functions constructed by
Wang \cite{W}:
 \[\sigma(t)=s_0(1+\delta-\frac \delta 2 e^{-\gamma t}),\ t\geq 0;\ \ Z(y)=\sin(\pi y),\ 0\leq y\leq 1,\]
 \[w(t,x)=Ke^{-\gamma t}Z(x/\sigma(t)),\ \ t\geq 0,\ 0\leq x\leq \sigma(t).\]
It is obvious that
  $$w(t,0)=w(t,\sigma(t))=0,\ \ \ \forall t>0.$$
Recall that $s_0<\pi\min\lk\{1,\sqrt{D/r}\rr\}$. Similar to the proof of Step 2 in Lemma 5.3 of \cite{W}, we can verify that, for the suitable small positive constants $\delta$ and
$\gamma$, and large positive constant $K$, the function $w$ satisfies
\[\left\{\begin{array}{lll}
 w_t-w_{xx}\geq w(1-w), &t>0, \ 0<x<\sigma(t),\\[1mm]
 w_t-Dw_{xx}\geq rw(1-w),&t>0, \ 0<x<\sigma(t),\\[1mm]
 w(0,x)\geq u_0(x), \ \ w(0,x)\geq v_0(x),\ & 0\leq x\leq s_0,
 \end{array}\right.\]
Moreover, for such fixed constants $\delta$, $\gamma$ and $K$, there exists $\mu_0>0$ such that
 $$\sigma'(t)+\mu(1+\rho)w_x(t,\sigma(t))\geq 0$$
for all $\mu\leq\mu_0$.

By Lemma 3.1, $s(t)\leq \sigma(t)$. Taking $t\to\infty$, we have $s_\infty\sigma(\infty)=s_0(1+\delta)<\infty$.\ \ \ $\Box$

\begin{theo}\lbl{th3.2} \ Suppose that $s_0<\Lambda$. For any one of problems {\rm (NFB)} and {\rm (DFB)}, there exist $\mu^*\geq\mu_*>0$, depending on $(u_0,v_0,s_0)$, such that $s_\infty\leq \Lambda$ if $\mu\leq\mu_*$, and $s_\infty=\infty$ if $\mu>\mu^*$.
\end{theo}

{\bf Proof}. \  Remember Lemmas \ref{le3.3} and \ref{le3.4}, the proof is similar to that of Theorem 5.4 in \cite{WZ}. We omit the details.\ \ \
\fbox{}

\section{Discussion}
\setcounter{equation}{0}

In this paper, we have examined a Lotka-Volterra type competition model with free boundary
$x=s(t)$ for both species, which describes the movement process through the free
boundary. We envision that the two species initially occupy the region $[0,s_0]$ and have a tendency to expand their territory together. Then we extend some results of \cite{DLin} and \cite{KY} for one species case and simplify the conditions in \cite{GW} for (NFB) with weak competition case. The dynamic behavior is discussed. Let $\Lambda=\frac\pi 2\min\lk\{1,\sqrt{D/r}\rr\}$ for the problem (NFB), and $\Lambda=\pi\min\lk\{1,\sqrt{D/r}\rr\}$
for the problem (DFB). It is proved that:

(i) If the size of initial habitat is not less than $\Lambda$, or it is less than $\Lambda$ but the moving parameter/coefficient $\mu$ of the free boundary is greater than $\mu^*$ (it depends on the initial data $u_0,v_0,s_0$), then $s_\infty=\infty$. Moreover,

(ia) To the problem (NFB), the dynamic behaviors are same as those of the ODE system (\ref{2.12}) for the cases: $0<k,h<1$, $0<k<1\leq h$, $0<h<1\leq k$;

(ib) To the problem (DFB), if $0<k,h<1$, then $u(t,x)$ and $v(t,x)$ satisfy
\[
\dd\liminf_{t\to\infty}u(t,x)\geq\ud u(x), \ \ \limsup_{t\to\infty}u(t,x)\leq\bar u(x),\ \
\dd\liminf_{t\to\infty}v(t,x)\geq\ud v(x), \ \ \limsup_{t\to\infty}v(t,x)\leq\bar v(x)
\]
uniformly in any compact subset of $[0,\infty)$, where $\ud u,\ \ud v,\ \bar u$ and $\bar v$ are positive functions given by
Theorem \ref{th2.5}.

(ii) While if the size of initial habitat is less than $\Lambda$ and the moving parameter/coefficient $\mu$ of the free boundary is less than $\mu_*$, then $\lim_{t\to\infty}s(t)<\Lambda$, and $\lim_{t\to\infty}\|u(t,\cdot)\|_{C([0,s(t)])}= \lim_{t\to\infty}\|v(t,\cdot)\|_{C([0,s(t)])}=0$. That is, the two species will disappear eventually.

The above conclusions not only provide the spreading-vanishing dichotomy and criteria governing
spreading and vanishing, but also provided the long time behavior of $(u, v)$ for both problems. If the size of initial habitat is small, and the moving parameter is small enough, it causes no population can survive eventually, while they can coexist if the size of habitat or the moving parameter is large enough, regardless of initial population size. This phenomenon suggests that the size of the initial habitat and the moving parameter are important to the survival for the two species.

\end{document}